\documentclass{amsart}[11pt]

\usepackage{fullpage}
\usepackage{amsmath}
\usepackage{amsfonts}
\usepackage{amssymb}

\newtheorem{lemma}{Lemma}[section]
\newtheorem{proposition}{Proposition}[section]
\newtheorem{theorem}{Theorem}[section]

\theoremstyle{remark}
\newtheorem{remark}{Remark}[section]
\theoremstyle{definition}

\title{Decay estimates for the one--dimensional wave equation with an
inverse power potential}

\begin{document}

\author{Roland Donninger}
\address{University of Chicago, Department of Mathematics,
5734 South University Avenue, Chicago, IL 60637, U.S.A.}
\email{donninger@uchicago.edu}
\thanks{The first author 
is an Erwin Schr\"odinger Fellow of the 
FWF (Austrian Science Fund) Project No. J2843 and 
he wants to thank Peter C. Aichelburg for his support and Piotr Bizo\'n for helpful discussions.}

\author{Wilhelm Schlag}
\address{University of Chicago, Department of Mathematics,
5734 South University Avenue, Chicago, IL 60637, U.S.A.}
\email{schlag@math.uchicago.edu}
\thanks{The second author was partly supported by the National
Science Foundation DMS-0617854 and by a Guggenheim fellowship.}

\begin{abstract}
We study the wave equation on the real line with a potential that falls off 
like $|x|^{-\alpha}$ for $|x| \to \infty$ where $2 < \alpha \leq 4$. We prove that the solution decays pointwise like $t^{-\alpha}$ as 
$t \to \infty$ provided that there are no resonances at zero energy and no bound states.
As an application we consider the $\ell=0$ Price Law for Schwarzschild black holes. This paper is part of our investigations into
decay of linear waves on a Schwarzschild background, see~\cite{donninger}, \cite{dss}. 
\end{abstract}

\maketitle

\section{Introduction}

There is an extensive literature in linear dispersive equations devoted to the study of decay for Schr\"odinger and wave equations with a potential, see for example~\cite{schlag3} for a survey. However, there seems to be little interaction with the physical community where the corresponding problem goes by the name of ``tails''.  Based largely on numerical evidence and nonrigorous arguments, physicists predict the decay of
solutions to wave equations on the line with potentials decaying like $|x|^{-\alpha}$ as $|x|\to\infty$, see for example~\cite{ching}, \cite{bizontail}. 
However, from a mathematical point of view this field is still largely open.

In the present paper we obtain decay estimates for the one--dimensional 
wave equation
\begin{equation}
\label{eq_1dwave}
\psi_{tt}-\psi_{xx}+V(x)\psi=0 
\end{equation}
with a potential of the form $V(x)=|x|^{-\alpha}
[c_\pm+O(x^{-\beta})]$ as $x \to \pm \infty$ where $2 < \alpha \leq 4$, 
$\beta=\frac{1}{2}(\alpha-2)^2$ and $c_\pm \in \mathbb{R}$.
We further assume that $V \in C^{[\alpha]+1}(\mathbb{R})$, 
the $O$--term satisfies
$|O^{(k)}(x)|\lesssim |x|^{-\beta-k}$ as $|x| \to \infty$  
for $k=0,1,\dots,[\alpha]+1$ and $V$ has
no bound states and no resonance at zero energy.
The symbol $[\alpha]$ denotes the smallest integer not less than
$\alpha$.
Under these assumptions we prove that the time evolution decays
pointwise like $t^{-\alpha}$ as $t \to \infty$ which confirms previous heuristics and numerics in the physical literature.
The precise statement is given in the theorem below in the form of weighted
$L^1$ to $L^\infty$ bounds for the sine and cosine evolutions.
A prominent physical application of our result is the problem of radial 
wave evolution in the presence of a Schwarzschild black hole in general relativity or other theories of gravity like Ho\v{r}ava--Lifshitz.

The motivation for our work is twofold: on the one hand, we need to incorporate a sharp decay estimate for spherical waves on Schwarzschild in our framework developed in~\cite{donninger} which only yielded $t^{{-2}}$ for vanishing angular momentum (and $t^{-2\ell-2}$ for angular momentum~$\ell$). 
This is an important building block for the proof of the sharp $t^{-3}$ decay for the wave equation on Schwarzschild without symmetry assumptions on the data which is established in the companion 
paper~\cite{dss}.
In addition,  we also wanted to confirm the predictions by physicists concerning the decay laws for inverse power potentials. 
However, the methods of this paper are not  able to cover
 the whole scale of exponents $\alpha$ that have been studied by physicists and 
it should be clear from our proof that more sophisticated techniques are required in order to 
obtain sharp results in greater generality.

\subsection{The main theorem}
In order to prove the result we construct the spectral measure of the 
associated
one--dimensional Schr\"odinger operator. 
More precisely, we define $Af:=-f''+Vf$ in $L^2(\mathbb{R})$ with domain
$\mathcal{D}(A):=H^2(\mathbb{R})$.
Thus, Eq.~(\ref{eq_1dwave}) can be written as
\begin{equation}
\label{eq_1dwaveop}
\psi_{tt}(t)+A\psi(t)=0
\end{equation}
where $\psi$ is now interpreted as a function of $t$ taking values in
$L^2(\mathbb{R})$.
The nonresonant condition on $V$ means that there does not exist a globally
bounded function $f \in L^\infty(\mathbb{R})$ with $Af=0$.
It is well--known that the operator $A$ is self--adjoint and, 
since there are no bound states, its spectrum is 
purely absolutely continuous and given by
$\sigma(A)=\sigma_{ac}(A)=[0,\infty)$ (cf. e.g. \cite{teschl}).
Thus, the functional calculus for self--adjoint operators yields 
the solution 
$$ \psi(t)=\cos(t\sqrt{A})f+
\frac{\sin(t\sqrt{A})}
{\sqrt{A}}g $$
of Eq.~(\ref{eq_1dwaveop}) with initial data
$(\psi(0),\psi_t(0))=(f,g)$.
The point is that the associated spectral measure can be expressed
via the Green's function $G(x,x',\lambda)$
which is the
kernel of the resolvent operator $((\lambda+i0)^2-A)^{-1}$.
Explicitly, we have
$$ \frac{\sin(t\sqrt{A})}{\sqrt{A}}f(x)=-\frac{2}{\pi}\lim_{N \to \infty}
\int_\mathbb{R} \int_{1/N}^N \sin(t\lambda)\mathrm{Im}[G(x,x',\lambda)]d\lambda
f(x')dx' $$
and a similar statement holds for the cosine evolution.
We refer the reader to \cite{donninger} and \cite{teschl} for a derivation of
these well--known facts.
The main result of the present paper is the following.

\begin{theorem}
\label{thm_main}
Let $V \in C^{[\alpha]+1}(\mathbb{R})$ with 
$V(x)=|x|^{-\alpha}[c_\pm+O(|x|^{-\beta})]$ as $x \to \pm \infty$ where
$2<\alpha \leq 4$, $\beta=\frac{1}{2}(\alpha-2)^2$, 
$c_\pm \in \mathbb{R}$ and $|O^{(k)}(|x|^{-\beta})|\lesssim |x|^{-\beta-k}$ for 
$k=1,2,\dots,[\alpha]+1$. Denote by $A$ the self--adjoint Schr\"odinger operator 
$Af:=-f''+Vf$ in $L^2(\mathbb{R})$ 
and assume that $A$ has no bound states and no resonance at zero
energy. 
Then the following decay bounds hold \footnote{Here and throughout this work the symbol 
$\langle x \rangle$ denotes a smooth
function that equals $|x|$ for $|x| \geq 2$ and satisfies $\langle x \rangle \geq 1$ for
all $x \in \mathbb{R}$.}:
$$ \|\langle \cdot \rangle^{-\alpha-1} 
\cos(t\sqrt{A})f\|_{L^\infty(\mathbb{R})}
\lesssim \langle t \rangle^{-\alpha} \left (\|\langle \cdot
\rangle^{\alpha+1}f'\|_{L^1(\mathbb{R})}+ 
\|\langle \cdot
\rangle^{\alpha+1}f\|_{L^1(\mathbb{R})}\right ) $$
and
$$ \left \|\langle \cdot \rangle^{-\alpha-1} \frac{\sin(t\sqrt{A})}{\sqrt{A}}
f \right \|_{L^\infty(\mathbb{R})}
\lesssim \langle t \rangle^{-\alpha} \|\langle \cdot
\rangle^{\alpha+1}f\|_{L^1(\mathbb{R})} $$
for all $t \geq 0$.
\end{theorem}

\begin{remark}
 As usual, if $V \in C^{[\alpha]+2}(\mathbb{R})$, one may improve the cosine estimate in Theorem \ref{thm_main} to
$$ \|\langle \cdot \rangle^{-\alpha-1} 
\cos(t\sqrt{A})f\|_{L^\infty(\mathbb{R})}
\lesssim \langle t \rangle^{-\alpha-1} \left (\|\langle \cdot
\rangle^{\alpha+1}f'\|_{L^1(\mathbb{R})}+ 
\|\langle \cdot
\rangle^{\alpha+1}f\|_{L^1(\mathbb{R})}\right ). $$
\end{remark}

\subsection{Application to the wave evolution on Schwarzschild}
The problem of radial wave evolution in the Schwarzschild geometry can be
reduced to Eq.~(\ref{eq_1dwave}) with the Regge--Wheeler potential
$$ V(x)=\frac{2M\sigma}{r^3(x)}\left (1-\frac{2M}{r(x)} \right ) $$
where $r(x)$ is given implicitly by
$$ x=r+2M\log \left ( \frac{r}{2M}-1 \right ), $$
the so--called tortoise coordinate, and $\sigma$ is a parameter that takes the
values $0,1,-3$ in different physical situations (however, since we restrict
ourselves to vanishing angular momentum in this paper, only $\sigma=1$ has physical meaning).
This is a prominent problem in classical general relativity.
We refer the reader to \cite{donninger} and references therein
for the history of waves on Schwarzschild and more background.
The asymptotics of the Regge--Wheeler potential are $V(x)=2M
x^{-3}+O(x^{-4}\log x)$ as $x \to \infty$ and, for $x \to -\infty$, $V(x)$ decays
exponentially.
Clearly, $V$ satisfies the requirements of Theorem \ref{thm_main} and therefore,
the solution decays pointwise like $t^{-3}$ as $t \to \infty$.
This is the famous $\ell=0$ Price Law for Schwarzschild black holes (cf.
\cite{Price1}, \cite{Price2}).
Note also that our estimates are sharp as far as the number of required
derivatives on the data is concerned.
While this paper was being written up,
Tataru~\cite{tataru} posted a preprint where the
sharp pointwise decay $t^{-3}$ for the wave equation on a 
very general asymptotically flat background is proved.
This remarkable result applies to both Schwarzschild and Kerr.
The best previously available result in spherical symmetry 
is due to Dafermos and Rodnianski
\cite{D-Rod} where $t^{-3+\varepsilon}$ decay has been proved.
We also mention the preprint \cite{Kron2} where the sharp $t^{-3}$ decay for
compactly supported data is
shown. However, the argument there is very involved and  the dependence of various constants on the initial data is unclear.

\subsection{Outline of the proof}
It is a common feature (see \cite{schlag3}) 
in dispersive estimates that the most important
contributions come from small energies. 
Therefore, one needs to develop a good understanding of the Green's function
$G(x,x',\lambda)$ around $\lambda=0$.
In other words, we have to obtain precise asymptotics of the Jost functions
$f_\pm$, defined by
$$ Af_\pm(\cdot,\lambda)=\lambda^2 f_\pm(\cdot,\lambda) $$
and $f_\pm(x,\lambda) \sim e^{\pm i\lambda x}$ for $x \to \pm \infty$, in the
limit $\lambda \to 0$.
These solutions are known to exist whenever $V \in L^1(\mathbb{R})$ (see
\cite{deift}), however, they are in general not smooth at $\lambda=0$ if the
potential has an inverse power decay.
It turns out that exactly this failure of smoothness 
is responsible for the characteristic decay 
$t^{-\alpha}$ given in Theorem~\ref{thm_main}.
In order to obtain the asymptotics of $f_\pm$ we use two different
representations of the Jost solutions.
First, we use the standard Volterra equation that defines the Jost solutions and
obtain an asymptotic expansion for $\lambda \to 0$ at 
$x=\pm \lambda^{-2/\alpha}$ (this is the turning point of the equation).
Second, we construct appropriate fundamental systems
of $Af=\lambda^2 f$ by perturbing in
$\lambda$ around $\lambda=0$.
We match the Jost solutions to these fundamental systems at $x=\pm
\lambda^{-2/\alpha}$
which allows us to obtain the precise asymptotics of the Wronskian
$W(f_-(\cdot,\lambda),f_+(\cdot,\lambda))$ in the limit $\lambda \to 0$.
At this point it is crucial that zero energy is not resonant, i.e., that
$W(f_-(\cdot,\lambda),f_+(\cdot,\lambda))$ does not vanish as $\lambda \to 0$.
Let $e(\lambda;x,x')$ denote the spectral measure of~$A$ at energy~$\lambda^{2}$. 
A basic conclusion of this paper is that
\begin{equation}\label{eq:spectralmeasure}
e(\lambda;0,0)=a\lambda + O(\lambda^{\alpha-1}) \qquad \lambda\to0
\end{equation}
where $a$ is some constant and the $O(\cdot)$ behaves like a symbol under differentiation in~$\lambda$. For the case of $\alpha=3$ this represents an improvement
over~\cite{donninger}  where only $e(\lambda;0,0)=O(\lambda)$ was shown. While the latter gave a decay
of $t^{-2}$ for the fundamental solution of the  wave equation (and nothing better), \eqref{eq:spectralmeasure} implies 
\[
\frac{\sin(t\sqrt{A})}{\sqrt{A}}(0,0) = 2\int_{0}^{\infty}  \frac{\sin(t\lambda)}{\lambda}   e(\lambda;0,0)\,\lambda d\lambda=O(t^{-3})
\]
as follows by means of three integrations by parts, see Lemma~\ref{lem_osc} below.
Note that the $a\lambda$ term in~\eqref{eq:spectralmeasure} does not contribute any unwanted boundary terms and simply drops
out at the third integration by parts. On the other hand, $\partial_{\lambda}^3 O(\lambda^{2})=O(\lambda^{-1})$ makes
a bounded contribution to the integral above, see the aforementioned Lemma~8.1. 
More generally, our methods yield an asymptotic representation of~$e(\lambda;x,x')$ for small~$\lambda$
and any $x,x'\in{\mathbb R}$ similar to~\eqref{eq:spectralmeasure}, which then implies the desired decay
bounds by means of oscillatory integral estimates  as in Section~8 of~\cite{donninger}, see also 
\cite{schlag1}, \cite{schlag2}. 

\subsection{Notations and conventions}
Throughout this work, $\alpha \in (2,4]$ is a fixed number and all implicit and explicit constants as well as all functions may depend on $\alpha$. However, we omit this dependence in the notation in order to improve readability.
The symbol $O(f(x))$ denotes a generic \emph{real--valued} function which is bounded by $|f(x)|$ in a domain of $x$ that follows from the context.
We write $O_\mathbb{C}(f(x))$ for complex--valued functions.
Furthermore, we say that $O(x^\gamma)$, $\gamma \in \mathbb{R}$, \emph{behaves like a symbol} if $O^{(k)}(x^\gamma)=O(x^{\gamma-k})$ for $k=1,2,\dots,[\alpha]+1$.
Finally, the letter $C$ (possibly with indices) denotes a positive constant that may change its value from line to line.

\section{Asymptotic expansions of the Jost functions}

We obtain asymptotic expansions of the Jost solutions $f_\pm(x,\lambda)$ at the turning point as
$\lambda \to 0$. As the potential has limited decay, this requires careful book keeping of various derivates.
In particular, we have found it advantageous to rescale the  usual Volterra equations by~$\lambda$
in order to control the derivatives with respect to~$\lambda$. 

\subsection{Expansion for $f_\pm(\pm \lambda^{-2/\alpha},\lambda)$}
The functions $f_\pm(\cdot,\lambda)$ are defined by 
$$ Af_\pm(\cdot,\lambda)=\lambda^2 f_\pm(\cdot,\lambda) $$
and the condition $f_\pm(x,\lambda) \sim e^{\pm i \lambda x}$ as $x \to \pm
\infty$.
By a straightforward application of the variation of constants formula we obtain
the Volterra integral equation
\begin{equation}
\label{eq_Jostint}
f_\pm(x,\lambda)=e^{\pm i \lambda x}+\int_x^{\pm
\infty}\frac{\sin(\lambda(y-x))}{\lambda}V(y)f_\pm(y,\lambda)dy.
\end{equation}
However, it turns out that it is more convenient to remove the oscillation and work
with the functions $m_\pm(x,\lambda):=e^{\mp i\lambda x}f_\pm(x,\lambda)$
instead.
The reason is that the functions $m_\pm$ behave well under differentiation as the following
lemma shows.

\begin{lemma}
\label{lem_symbf+}
The functions $m_\pm(x,\lambda):=e^{\mp i\lambda x}f_\pm(x,\lambda)$ satisfy
the estimates
$$ |\partial_\lambda^k \partial_x^j \left [m_\pm(x,\lambda)-1 \right ]|
\lesssim \langle x
\rangle^{-(\alpha-2)-j}\lambda^{-k} $$
for $\pm x \geq 0$, small $\lambda>0$ and $0 \leq j+k \leq [\alpha]+1$.
\end{lemma} 

\begin{proof}
The function $m_+$ satisfies the Volterra equation
\begin{align*}
m_+(x,\lambda)&=1+\frac{1}{2i\lambda}\int_x^\infty \left (e^{2i\lambda(y-x)}-1
\right )V(y)m_+(y,\lambda)dy \\
&=1+\frac{1}{2i}\int_0^\infty \left (e^{2i\eta}-1 \right)\lambda^{-2}V(\eta
\lambda^{-1}+x)m_+(\eta \lambda^{-1}+x,\lambda)d\eta
\end{align*}
as follows directly from Eq.~(\ref{eq_Jostint}).
We set $M:=\{0,1,\dots,[\alpha]+1\} \times \{0,1,\dots,[\alpha]+1\}$ and 
define a bijection $n: 
M \to \{0,1,\dots,|M|-1 \}$
according to (for e.g. $[\alpha]=3$)
$$ \left [\begin{array}{ccccc}
n(0,0) & n(1,0) & n(2,0) & n(3,0) & n(4,0) \\
n(0,1) & n(1,1) & n(2,1) & n(3,1) & \\
n(0,2) & n(1,2) & n(2,2) & &  \\
n(0,3) & n(1,3) &  & & \\
n(0,4) & & & & 
\end{array} \right ]= 
\left [\begin{array}{ccccc}
0 & 1 & 3 & 6 & 10 \\ 
2 & 4 & 7 & 11 &\\
5 & 8 & 12 & & \\
9 & 13 & & & \\
14 & & & &
\end{array} \right ]. $$
The existence theorem for Volterra equations 
(see e.g.~\cite{deift}, \cite{schlag1}, \cite{schlag2} or \cite{donninger}) shows that the lemma is true for
$(j,k)=(0,0)$. Fix $(j,k) \in M$ with $n(j,k) \leq |M|-2$ 
and assume the lemma is true for all $(\ell,m)
\in M$ with $n(\ell,m)\leq n(j,k)$. 
We show that this implies the claim for $(j',k') \in M$ where
$n(j',k')=n(j,k)+1$.
There are two possibilities: either $(j',k')=(k+1,0)$ (if $j=0$) or
$(j',k')=(j-1,k+1)$.
In the former case we have
\begin{align*} \partial_x^{k+1}m_+(x,\lambda)&=\frac{1}{2i\lambda^2}
\sum_{\ell=0}^{k+1} 
\left (\begin{array}{c}k+1 \\ \ell \end{array} \right )
\int_0^\infty \left (e^{2i\eta}-1 \right )\partial_x^{k+1-\ell}V(\eta
\lambda^{-1}+x)\partial_x^\ell m_+(\eta \lambda^{-1}+x,\lambda)d\eta \\
&=O_\mathbb{C}(\langle x \rangle^{-(\alpha-2)-(k+1)})+
\frac{1}{2i\lambda}\int_x^\infty \left (e^{2i\lambda(y-x)}-1
\right )V(y)\partial_y^{k+1}m_+(y,\lambda)dy 
\end{align*}
by assumption and the properties of $V$ (recall that 
$|m_+(x,\lambda)|\lesssim 1$).
Thus, the standard result on Volterra equations yields
$|\partial_x^{k+1} m_+(x,\lambda)| \lesssim \langle x \rangle^{-(\alpha-2)-(k+1)}$ for
all $x \geq 0$.
For the second case it is useful to note that
$$ \left |\partial_\lambda^m \partial_x^\ell m_+(\eta
\lambda^{-1}+x,\lambda)\right |
\lesssim \langle \eta \lambda^{-1}+x \rangle^{-(\alpha-2)-\ell} \lambda^{-m} $$
if $n(\ell,m) \leq n(j,k)$, $(\ell,m) \not= (0,0)$ by the chain rule (cf. also
\cite{donninger}, Lemma~9.1).
This and the properties of $V$ imply that
\begin{align*}
\partial_\lambda^{k+1} \partial_x^{j-1} m_+(x,\lambda)&=\sum_{m=0}^{k+1} \left (
\begin{array}{c}k+1 \\ m \end{array} \right )\int_0^\infty \left (e^{2i\eta}-1
\right ) \partial_x^{j-1} \left [ \partial_\lambda^{k+1-m}
\frac{V(\eta \lambda^{-1}+x)}{2i\lambda^2}\partial_\lambda^m m_+(\eta
\lambda^{-1}+x,\lambda) \right ]d\eta \\
&=O_\mathbb{C}(\langle x \rangle^{-(\alpha-2)-(j-1)}\lambda^{-(k+1)})+
\frac{1}{2i\lambda}\int_x^\infty \left (e^{2i\lambda(y-x)}-1
\right )V(y)\partial_\lambda^{k+1}\partial_y^{j-1}m_+(y,\lambda)dy 
\end{align*}
and we obtain $|\partial_\lambda^{k+1} \partial_x^{j-1} m_+(x,\lambda)| \lesssim
\langle x \rangle^{-(\alpha-2)-(j-1)} \lambda^{-(k+1)}$ as claimed.
The proof for $m_-$ is identical.
\end{proof}

\begin{lemma}
\label{lem_f+}
The Jost solutions have the asymptotic expansion
\begin{align*} f_\pm(\pm
\lambda^{-2/\alpha},\lambda)&=
1+i \mu+\frac{1}{2}(c_\pm-1)\mu^2
-ic_\pm \mu^3 \log \mu+O_\mathbb{C}(\mu^3) \mbox{ if }\alpha=3 \\
 f_\pm(\pm
\lambda^{-2/\alpha},\lambda)&=
1+i\mu+\left (\frac{c_\pm}{(\alpha-1)(\alpha-2)}-\frac{1}{2} \right )\mu^2 + i
\left (\frac{c_\pm}{(\alpha-2)(\alpha-3)}-\frac{1}{6}\right )\mu^3+O_\mathbb{C}(\mu^\alpha) \mbox{ if } \alpha \not=3 
\end{align*}
for small $\lambda>0$ where $\mu=\lambda^{1-2/\alpha}$ and the $O$--terms behave like symbols
under differentiation.
\end{lemma}

\begin{proof}
We only prove the assertion for $f_+$ since the rest follows by means of symmetry
considerations.
Again, we work with $m_\pm$ and rewrite Eq.~(\ref{eq_Jostint}) as
\begin{align}
\label{eq_proofJost}
m_+(x,\lambda)-1&=\frac{1}{2i\lambda}\int_x^\infty
\left ( e^{2i\lambda(y-x)}-1 \right )V(y)dy \\
&+\underbrace{\frac{1}{2i\lambda}\int_x^\infty
\left (e^{2i\lambda(y-x)}-1 \right )V(y)\left [m_+(y,\lambda)-1 \right] dy}_{
=:A(x,\lambda)} \nonumber
\end{align}
and the existence theorem for Volterra equations yields 
$ |m_+(x,\lambda)-1|\lesssim \langle x \rangle^{-(\alpha-2)} $
for all $x \geq 0$ since
$$ \left | \frac{1}{2i\lambda}\int_x^\infty
\left ( e^{2i\lambda(y-x)}-1 \right )V(y)dy \right |\lesssim
\int_x^\infty \langle y \rangle^{-\alpha+1}dy \lesssim \langle x \rangle^{-\alpha+2}. $$
This implies $A(x,\lambda)=O_\mathbb{C}(\langle x \rangle^{-2(\alpha-2)})$
and we conclude that
\begin{equation}
\label{eq_Jostproofm+}
m_+(x,\lambda)=1+\frac{1}{2i\lambda}\int_x^\infty
\left ( e^{2i\lambda(y-x)}-1 \right )V(y)dy+O_\mathbb{C}(\langle x
\rangle^{-2(\alpha-2)}). 
\end{equation}
Furthermore, for $0 \leq j+k \leq [\alpha]+1$, we have
\begin{align*} \partial_\lambda^k \partial_x^j A(x,\lambda)&=
\partial_\lambda^k \partial_x^j \frac{1}{2i}\int_0^\infty
\left (e^{2i\eta}-1 \right )\lambda^{-2}V(\eta \lambda^{-1}+x)
\left [m_+(\eta \lambda^{-1}+x,\lambda)-1 \right] d\eta \\
&=\int_0^\infty \left (e^{2i\eta}-1 \right )
O_\mathbb{C}(\langle \eta \lambda^{-1}+x \rangle^{-2\alpha+2-j}\lambda^{-2-k})d\eta \\
&=\lambda^{-1-k}\int_x^\infty \left (e^{2i\lambda(y-x)}-1 \right)O_\mathbb{C}
(\langle y \rangle^{-2\alpha+2-j})dy=O_\mathbb{C}(\langle x \rangle^{-2(\alpha-2)-j}\lambda^{-k})
\end{align*}
by Lemma \ref{lem_symbf+} and thus, the $O_\mathbb{C}$--term in
Eq.~(\ref{eq_Jostproofm+}) behaves like a symbol with respect to differentiation
in $x$ and $\lambda$.
By using $V(x)=c_+ x^{-\alpha}+O(x^{-\alpha-\frac{1}{2}(\alpha-2)^2})$ for $x \geq 2$, 
Eq.~(\ref{eq_Jostproofm+}) reduces to
\begin{align}
\label{eq_Jostproofm+2}
m_+(x,\lambda)=1+\frac{c_+}{2i\lambda}\int_x^\infty
\left ( e^{2i\lambda(y-x)}-1 \right )y^{-\alpha}dy+O_\mathbb{C}(x^{-\frac{1}{2}\alpha(\alpha-2)})
+O_\mathbb{C}(x^{-2(\alpha-2)})
\end{align}
and again, by the properties of $V$ and the same argument as above for
$A(x,\lambda)$, we obtain the symbol behavior of the $O$--terms in
Eq.~(\ref{eq_Jostproofm+2}).

Specializing to $\alpha=3$, evaluation at $x=\lambda^{-2/3}$ yields
\begin{equation}
\label{eq_Jostproofeval} 
m_+(\lambda^{-2/3},\lambda)=1+\frac{c_+}{2i\lambda}\int_{\lambda^{-2/3}}^\infty
\left ( e^{-2i\lambda^{1/3}}e^{2i\lambda y}-1 \right )y^{-3}dy+O_\mathbb{C}(\lambda)
\end{equation}
where the $O$--term behaves like a symbol under differentiation.
Now we have
\begin{align*} 
\frac{1}{2i\lambda}e^{-2i\lambda^{1/3}}
\int_{\lambda^{-2/3}}^\infty e^{2i\lambda y}y^{-3}dy&=
\frac{\lambda}{2i}e^{-2i\lambda^{1/3}}\int_{\lambda^{1/3}}^\infty
e^{2i\eta}\eta^{-3}d\eta\\
&=\frac{\lambda}{2i}e^{-2i\lambda^{1/3}}\left (\int_{\lambda^{1/3}}^1
e^{2i\eta}\eta^{-3}d\eta+\int_1^\infty e^{2i\eta}\eta^{-3}d\eta \right ) \\
&=\frac{\lambda}{2i}e^{-2i\lambda^{1/3}}\left (\sum_{j=0}^\infty
\frac{(2i)^j}{j!}\int_{\lambda^{1/3}}^1 
\eta^{j-3}d\eta+c_1 \right ) \\
&=-\frac{1}{4}i\lambda^{1/3}+\frac{1}{2}\lambda^{2/3}-\frac{1}{3}i\lambda \log
\lambda+O_\mathbb{C}(\lambda)
\end{align*}
and clearly, the $O$--terms behave like
symbols under differentiation since they stem from the Taylor series expansion
of the exponential.
Inserting this in Eq.~(\ref{eq_Jostproofm+2}) we obtain
$$ m_+(\lambda^{-2/3},\lambda)=1+\frac{c_+}{2}\lambda^{2/3}
-\frac{c_+}{3}i\lambda \log
\lambda+O_\mathbb{C}(\lambda) $$
and expanding
$f_+(\lambda^{-2/3},\lambda)=e^{i\lambda^{1/3}}m_+(\lambda^{-2/3},\lambda)$
finishes the proof for $\alpha=3$.

For $\alpha \not=3$ we evaluate Eq.~(\ref{eq_Jostproofm+2}) at $x=\lambda^{-2/\alpha}$ which yields
\begin{equation*}
m_+(\lambda^{-2/\alpha},\lambda)=1+\frac{c_+}{2i\lambda}\int_{\lambda^{-2/\alpha}}^\infty
\left ( e^{-2i\lambda^{1-2/\alpha}}e^{2i\lambda y}-1 \right )y^{-\alpha}dy+O_\mathbb{C}(\lambda^{\alpha(1-2/\alpha)})+O_\mathbb{C}(\lambda^{4(1-2/\alpha)})
\end{equation*}
and performing analogous Taylor series expansions as in the case $\alpha=3$ yields the claim.
\end{proof}

\subsection{Expansion for $f_\pm'(\pm \lambda^{-2/\alpha},\lambda)$}
In order to calculate Wronskians we also need expansions for the derivatives of
the Jost functions.

\begin{lemma}
\label{lem_f+der}
The functions $f_\pm'(\pm \lambda^{-2/\alpha},\lambda)$ have the asymptotic expansions
$$ f_\pm'(\pm \lambda^{-2/\alpha},\lambda)=\pm \lambda^{2/\alpha}\left
[i\mu-\left (\frac{c_\pm}{\alpha-1}+1 \right )\mu^2-i\left
(\frac{c_\pm}{\alpha-2}+\frac{1}{2} \right )\mu^3+O(\mu^\alpha)+i O(\mu^{\alpha+1}) \right ] $$
for small $\lambda>0$ where $\mu=\lambda^{1-2/\alpha}$ and the $O$--terms behave like symbols.
\end{lemma}

\begin{proof}
We use the representation
\begin{align*} m_+'(x,\lambda)=&-c_+\int_x^\infty e^{2i\lambda(y-x)}y^{-\alpha}dy 
-\int_x^\infty e^{2i\lambda(y-x)}O(y^{-\alpha-\frac{1}{2}(\alpha-2)^2})dy\\
&-\underbrace{\int_x^\infty e^{2i\lambda(y-x)}
V(y)\left [m_+(y,\lambda)-1 \right ] dy}_{=:B(x,\lambda)}
\end{align*}
which follows directly from Eq.~(\ref{eq_proofJost}) by differentiation.
Lemma \ref{lem_symbf+} shows that 
$$ B(x,\lambda)=O_\mathbb{C}(\langle x \rangle^{-2(\alpha-2)-1}) $$
 and the $O$--term behaves like a symbol.
Furthermore, for $\alpha \in (3,4]$ one may refine this bound to
$$ B(x,\lambda)=O(\langle x \rangle^{-2(\alpha-2)-1})+iO(\langle x \rangle^{-2(\alpha-2)}\lambda) $$
by noting that
\begin{align*}\mathrm{Im}\left [ m_+(x,\lambda) \right ]&=-\frac{1}{2\lambda}\int_x^\infty \left [\cos(2\lambda(y-x))-1 \right ]V(y)\mathrm{Re} \left [ m_+(y,\lambda) \right ]dy \\
&+\frac{1}{2\lambda}\int_x^\infty  \sin(2\lambda(y-x))V(y)\mathrm{Im}\left [m_+(y,\lambda) \right ]dy \\
&=O(\langle x \rangle^{-(\alpha-3)}\lambda)+\frac{1}{2\lambda}\int_x^\infty  \sin(2\lambda(y-x))V(y)\mathrm{Im} \left [m_+(y,\lambda) \right ] dy
\end{align*}
which implies $\mathrm{Im} \left [m_+(x,\lambda) \right ]=
O(\langle x \rangle^{-(\alpha-3)}\lambda)$ by the Volterra existence theorem.
Thus,
$$ m_+'(\pm \lambda^{-2/\alpha},\lambda)=-c_+\int_{\lambda^{-2/\alpha}}^\infty e^{2i\lambda(y-\lambda^{-2/\alpha})}y^{-\alpha}dy 
+O(\lambda^{\alpha(1-2/\alpha)+2/\alpha})+iO(\lambda^{(\alpha+1)(1-2\alpha)+2/\alpha})
$$
where the $O$--terms behave like symbols and the claim follows by expanding
$$ f'(\lambda^{-2/\alpha},\lambda)=e^{i\lambda^{1-2/\alpha}}\left [i\lambda 
m_+(\lambda^{-2/\alpha},\lambda)+m'(\lambda^{-2/\alpha},\lambda) \right ] $$
as in the proof of Lemma \ref{lem_f+}.
\end{proof}

\section{Perturbation in energy}
In this section we construct two fundamental systems of solutions to
$Af=\lambda^2 f$ by perturbing in $\lambda$ around $\lambda=0$.
This is done in such a way that these fundamental systems can be matched to the
Jost functions.

\subsection{Zero energy solutions}
As a first step we obtain asymptotics for solutions of $Af=0$.

\begin{lemma}
\label{lem_zeroenergy}
There exists a fundamental system $\{u_0^\pm, u_1^\pm\}$ 
for the equation $Af=0$ with $W(u_0^\pm, u_1^\pm)=\mp 1$ and
\begin{align*} 
u_1^\pm(x)&=1+\frac{c_\pm}{(\alpha-1)(\alpha-2)}\langle x \rangle^{-(\alpha-2)}
+O(\langle x \rangle^{-\frac{1}{2}\alpha(\alpha-2)}) \\
u_0^\pm(x)&=\langle x \rangle \left [1-c_\pm\langle x \rangle^{-1}\log \langle x
\rangle+O(\langle x \rangle^{-3/2}) \right ] \mbox{ if }\alpha=3 \\
u_0^\pm(x)&=\langle x \rangle \left [1+\frac{c_\pm}{(\alpha-2)(\alpha-3)}\langle x \rangle^{-(\alpha-2)}
+O(\langle x \rangle^{-\frac{1}{2}\alpha(\alpha-2)}) \right ] \mbox{ if }\alpha\not=3 \\
\end{align*}
for $\pm x \geq 0$ where the $O$--terms behave like symbols under
differentiation.
\end{lemma}

\begin{proof}
As always we only prove the $+$ case. We consider the equation
$$ u_1^+(x)=1+\int_x^\infty(y-x)V(y)u_1^+(y)dy $$
and the Volterra existence theorem shows that there is a unique solution
satisfying $|u_1^+(x)|\lesssim 1$ and thus, the equation implies
$|u_1^+(x)-1|\lesssim \langle x \rangle^{-(\alpha-2)}$ since $|V(y)|\lesssim \langle y
\rangle^{-\alpha}$.
By a straightforward calculation it follows that $Au_1^+=0$.
For $x \geq 2$ we rewrite the Volterra equation for $u_1^+$ as 
\begin{align*}
u_1^+(x)-1&=c_+ \int_x^\infty(y-x)y^{-\alpha} dy+\int_x^ \infty(y-x)O(y^{-\alpha-\frac{1}{2}\alpha(\alpha-2)^2})dy \\
&+\int_x^\infty(y-x)V(y)\left [u_1^+(y)-1 \right
]dy \\
&=\frac{c_+}{(\alpha-1)(\alpha-2)}x^{-(\alpha-2)}+O(\langle x \rangle^{-\frac{1}{2}\alpha(\alpha-2)})
\end{align*}
since $\langle x \rangle^{-2(\alpha-2)} \lesssim \langle x \rangle^{-\frac{1}{2}\alpha(\alpha-2)}$.
The $O$--term behaves like a symbol due to the properties of $V$ and a
simple induction.
Furthermore, the solution can be extended to $x \geq 0$ since $V \in C^{[\alpha]+1}(\mathbb{R})$.

In order to construct $u_0^+$ in the domain $x \geq x_1$, 
we use the reduction ansatz
$$ u_0^+(x)=u_1^+(x)\left (\int_{x_1}^x u_1^+(y)^{-2}dy+\alpha_+
\right ) $$
where $x_1$ is chosen so large that 
$|u_1^+(x)| > 0$ for
$x \geq x_1$ and $\alpha_+$ is a real constant.
Inserting the asymptotics of $u_1^+$ and choosing $\alpha_+$ appropriately,
we obtain $u_0^+$ with the claimed asymptotics for $x \geq x_1$ and, by extension, for $x \geq 0$.
The Wronskian condition is satisfied by construction.
\end{proof}

\subsection{Perturbation in energy}
Now we perturb in energy, i.e., we treat the right--hand side of $Af=\lambda^2
f$ as a perturbation.

\begin{lemma}
\label{lem_pertenergy}
The equation $Af=\lambda^2 f$ has fundamental systems
$\{u_0^\pm(\cdot,\lambda),u_1^\pm(\cdot,\lambda)\}$ with
$W(u_0(\cdot,\lambda),u_1(\cdot,\lambda))=\mp 1$ and
\begin{align*}
u_0^\pm(x,\lambda)&=u_0^\pm(x)-\frac{1}{6}\langle x \rangle^3 \lambda^2 
+O(\lambda^2 \langle x \rangle^2 \log \langle x \rangle)+O(\langle x \rangle^5
\lambda^4) \\
u_1^\pm(x,\lambda)&=u_1^\pm(x)-\frac{1}{2}\langle
x \rangle^2 \lambda^2+O(\lambda^2 \langle x
\rangle \log \langle x \rangle)+O(\langle x \rangle^4 \lambda^4) 
\end{align*}
if $\alpha=3$ as well as
\begin{align*}
u_0^\pm(x,\lambda)&=u_0^\pm(x)-\frac{1}{6}\langle x \rangle^3 \lambda^2 
+O(\langle x \rangle^{5-\alpha}\lambda^2)+O(\langle x \rangle^5
\lambda^4) \\
u_1^\pm(x,\lambda)&=u_1^\pm(x)-\frac{1}{2}\langle
x \rangle^2 \lambda^2+O(\langle x
\rangle^{4-\alpha}\lambda^2)+O(\langle x \rangle^4 \lambda^4) 
\end{align*}
if $\alpha \not=3$.
These expansions are valid for small $\lambda>0$, $x \in [0,\lambda^{-1}]$ (resp. $x \in
[-\lambda^{-1},0]$) and the $O$--terms behave like symbols
under differentiation with respect to $x$ and $\lambda$.
\end{lemma}

\begin{proof}
By symmetry considerations it suffices to prove the $+$ case. 
We define an integral kernel 
$$ K_0(x,y):=\langle y \rangle u_0^+(y)\frac{u_1^+(x)}{\langle
x\rangle}-\frac{u_0^+(x)}{\langle x \rangle}\langle y \rangle u_1^+(y) $$
for $x,y\geq 0$.
A straightforward calculation shows that the solution $u_0^+(\cdot,\lambda)$ of
\begin{equation}
\label{eq_proofpertenergyu0}
\frac{u_0^+(x,\lambda)}{\langle x \rangle}=\frac{u_0^+(x)}{\langle
x\rangle}+\lambda^2 \int_0^x K_0(x,y)\frac{u_0^+(y,\lambda)}{\langle y
\rangle}dy 
\end{equation}
satisfies $Au_0^+(\cdot,\lambda)=\lambda^2 u_0^+(\cdot,\lambda)$.
Let $\lambda \in (0,\lambda_0)$ for some fixed $\lambda_0>0$.
By Lemma~\ref{lem_zeroenergy} we have the bound
$|K_0(x,y)|\lesssim \langle y \rangle$ for $0 \leq y \leq x$ and thus,
$$ \lambda^2 \int_0^{\lambda^{-1}} \sup_{x: 0 \leq y \leq x}|K_0(x,y)|dy \lesssim
1 $$
and the existence theorem for Volterra equations implies that
Eq.~(\ref{eq_proofpertenergyu0}) has a unique solution satisfying
$$ \left |\frac{u_0(x,\lambda)}{\langle x \rangle} \right |\lesssim 1 $$
for all $x \in [0,\lambda^{-1}]$.
Thus, Eq.~(\ref{eq_proofpertenergyu0}) shows that
$$ \left |\frac{u_0(x,\lambda)}{\langle x \rangle}-\frac{u_0(x)}{\langle x
\rangle}\right |=\left |
\lambda^2 \int_0^x K_0(x,y)\frac{u_0^+(y,\lambda)}{\langle y
\rangle}dy \right | \lesssim \langle x \rangle^2 \lambda^2 $$
and we conclude
\begin{align}
\label{eq_proofpertenergyu02}
\frac{u_0(x,\lambda)}{\langle x \rangle}-\frac{u_0(x)}{\langle x
\rangle}=&\lambda^2 \int_0^x K_0(x,y)\frac{u_0^+(y)}{\langle y
\rangle}dy \\
&+\underbrace{\lambda^2 \int_0^x K_0(x,y)\left [ 
\frac{u_0^+(y, \lambda)}{\langle y \rangle}-
\frac{u_0^+(y)}{\langle y \rangle} \right ]dy}_{=O(\langle x \rangle^4
\lambda^4)} \nonumber
\end{align}
and by an appropriate induction (cf.~the proof of Lemma~\ref{lem_f+} or
Proposition~4.1 in~\cite{donninger}) it follows that the $O$--term behaves like
a symbol with respect to differentiation in $x$ and $\lambda$.
By using the expansions from Lemma~\ref{lem_zeroenergy} we
obtain
$$ \lambda^2 \int_0^x K_0(x,y)\frac{u_0^+(y)}{\langle y
\rangle}dy=-\frac{1}{6}\langle x \rangle^2\lambda^2 + O(\langle x
\rangle^{4-\alpha}\lambda^2) $$
if $\alpha \not= 3$ and similarly,
$$ \lambda^2 \int_0^x K_0(x,y)\frac{u_0^+(y)}{\langle y
\rangle}dy=-\frac{1}{6}\langle x \rangle^2\lambda^2 + O(\lambda^2 \langle x \rangle 
\log \langle x \rangle) $$
in the case $\alpha=3$. Again, the $O$--terms inherit the symbol behavior with respect to differentiation in both $x$ and $\lambda$.
Inserting this into Eq.~(\ref{eq_proofpertenergyu02}) yields the claim for
$u_0^+(\cdot,\lambda)$.

For $u_1^+$ we proceed along the same lines and use the Volterra equation
$$ u_1^+(x,\lambda)=u_1^+(x)+\lambda^2 \int_0^x
K_1(x,y)u_1^+(y,\lambda)dy $$
where now 
$$ K_1(x, y):=u_0^+(y)u_1^+(x)-u_0^+(x)u_1^+(y) $$
for $x,y \geq 0$.
Repeating the above arguments for this equation yields the claim.
The Wronskian condition follows from evaluation at $x=0$.
\end{proof}	

\subsection{Evaluation at $x=\pm \lambda^{-2/\alpha}$}
In order to match the  fundamental system $\{u_0^\pm(\cdot,\lambda),
u_1^\pm(\cdot,\lambda)\}$ to the Jost functions $f_\pm(\cdot,\lambda)$ one needs
to evaluate the functions $u_j^\pm(x,\lambda)$ at $x=\pm \lambda^{-2/\alpha}$.

\begin{lemma}
\label{lem_pertenergyeval}
The functions $u_j^\pm(\pm \lambda^{-2/\alpha},\lambda)$, $j=0,1$, have the
asymptotic expansions
\begin{align*}
u_0^\pm(\pm \lambda^{-2/3},\lambda)&=\lambda^{-2/3} \left [1+2 c_\pm
\mu^2\log \mu-\frac{1}{6}\mu^2+O(\mu^3) \right ] \mbox{ if }\alpha=3 \\
u_0^\pm(\pm \lambda^{-2/\alpha},\lambda)&=\lambda^{-2/\alpha}\left [1+\left
 (\frac{c_\pm}{(\alpha-2)(\alpha-3)}-\frac{1}{6} \right )\mu^2+O(\mu^\alpha) \right ] \mbox{ if } \alpha 
\not=3 \\
\partial_x u_0^\pm(\pm \lambda^{-2/\alpha},\lambda)&=\pm \left [1-\left (\frac{c_\pm}{\alpha-2}+\frac{1}{2} \right )\mu^2+O(\mu^\alpha) \right ] \\
u_1^\pm(\pm \lambda^{-2/\alpha},\lambda)&=1+\left (\frac{c_\pm}{(\alpha-1)(\alpha-2)}-\frac{1}{2} \right ) \mu^2+O(\mu^\alpha) \\
\partial_x u_1^\pm(\pm \lambda^{-2/\alpha},\lambda)&=\mp \lambda^{2/\alpha}\left [\left (\frac{c_\pm}{\alpha-1} +1 \right )\mu^2+O(\mu^\alpha) \right ]  
\end{align*}
where $\mu=\lambda^{1-2\alpha}$.
These expansions are valid for 
$\lambda>0$ small and all $O$--terms behave like symbols.
\end{lemma}

\begin{proof}
Based on Lemma~\ref{lem_pertenergy} this is just a straightforward computation.
Note that, thanks to the symbol behavior of the $O$--terms in Lemma 
\ref{lem_pertenergy}, the asymptotics for the 
$x$--derivatives can be obtained by formal
differentiation and evaluation at $x=\pm \lambda^{-2/\alpha}$.
\end{proof}

\section{Construction of the spectral measure at zero energy}

With the asymptotic expansions from the last section at hand we can now
calculate the Wronskian of the Jost solutions $f_-(\cdot,\lambda)$ and
$f_+(\cdot,\lambda)$. 

\begin{lemma}
\label{lem_connection}
Let $a_j^\pm(\lambda):=W(f_\pm(\cdot,\lambda),u_j^\pm(\cdot,\lambda))$.
Then we have the asymptotics
\begin{align*}
a_0^\pm(\lambda)&=\pm \left [1-ic_\pm \lambda \log
\lambda+O_\mathbb{C}(\lambda) \right ] \mbox{ if } \alpha=3\\
a_0^\pm(\lambda)&=\pm \left [ 1+O_\mathbb{C}(\mu^\alpha) \right ] 
\mbox{ if }\alpha \not=3 \\
a_1^\pm(\lambda)&=\pm \lambda^{2/\alpha} \left [-i\mu+O(\mu^\alpha)+iO(\mu^{\alpha+1}) \right ]
\end{align*}
for small $\lambda>0$  
where $\mu=\lambda^{1-2/\alpha}$ and all $O$--terms behave
like symbols.
\end{lemma}

\begin{proof}
By using the expansions from 
Lemmas \ref{lem_f+}, \ref{lem_f+der} and \ref{lem_pertenergy} 
this is a straightforward computation.
\end{proof}

Note that we have the representation
\begin{equation}
\label{eq_Jostrep}
f_\pm(x,\lambda)=\mp a_1^\pm(\lambda)u_0^\pm(x,\lambda)\pm
a_0^\pm(\lambda)u_1^\pm(x,\lambda) 
\end{equation}
for the Jost functions.
Next, we define the connection coefficients 
$$ b_{jk}:=W(u_j^-(\cdot,\lambda), u_k^+(\cdot,\lambda)) $$
and emphasize that $b_{jk}$ are real and 
independent of $\lambda$ (evaluate the Wronskians
at $x=0$ and recall the construction of $u_j^\pm(\cdot,\lambda)$ in Lemma~\ref{lem_pertenergy}).
Moreover, the nonresonant condition is equivalent to $b_{11} \not=0$.
In order to calculate the spectral measure we need to understand expressions
of the form
$$ \mathrm{Im}\left [\frac{f_-(x,\lambda)f_+(x',\lambda)}
{W(f_-(\cdot,\lambda),f_+(\cdot,\lambda))} \right ] $$
for small $\lambda$ where we write $W(f_-(\cdot,\lambda),f_+(\cdot,\lambda))=W(f_-,f_+)(\lambda)$ from now on.
In view of Eq.~(\ref{eq_Jostrep}) and since  $u_j^\pm(x,\lambda)$  are real-valued, 
we therefore have to study 
$$ 
\mathrm{Im}\left [ \frac{a_j^-(\lambda)a_k^+(\lambda)}{W(f_-,f_+)(\lambda)}
\right ]
$$
for small $\lambda$. The reader should compare the following proposition with~\eqref{eq:spectralmeasure}. 

\begin{proposition}
\label{prop_spectralmeasure}
There exists  $\lambda_0>0$ such that
$$ \mathrm{Im}\left [ \frac{a_j^-(\lambda)a_k^+(\lambda)}{W(f_-,f_+)(\lambda)}
\right ]=\alpha_{jk}\lambda+O(\lambda^{\alpha-1}) $$
for all $\lambda \in [0,\lambda_0]$ where $\alpha_{jk} \in \mathbb{R}$ and the
$O$--term behaves like a symbol.
\end{proposition} 

\begin{proof}
The Wronskian is given by
$$ W(f_-,f_+)(\lambda)=-b_{00}
a_1^-(\lambda)a_1^+(\lambda)+b_{10}a_0^-(\lambda)a_1^+(\lambda)
+b_{01}a_1^-(\lambda)a_0^+(\lambda)-b_{11}a_0^-(\lambda)a_0^+(\lambda), $$
see Eq.~(\ref{eq_Jostrep}) above.
We first consider the case $\alpha \not=3$ and from the expansions in Lemma \ref{lem_connection} we obtain
\begin{align*}
a_0^- a_0^+(\lambda)&=-1+O_\mathbb{C}(\mu^\alpha) \\
a_0^- a_1^+(\lambda)&=\lambda^{2/\alpha}\left [i\mu+O(\mu^\alpha)+iO(\mu^{\alpha+1}) \right ]\\
a_1^- a_0^+(\lambda)&=\lambda^{2/\alpha}\left [i\mu+O(\mu^\alpha)+iO(\mu^{\alpha+1}) \right ]\\
a_1^- a_1^+(\lambda)&=\lambda^{1+2/\alpha}\left [\mu+O(\mu^\alpha)+iO(\mu^{\alpha+1})) \right ]
\end{align*}
where, as before, $\mu=\lambda^{1-2/\alpha}$.
This implies
\begin{align*}
 \mathrm{Im}\left [\frac{a_0^- a_0^+(\lambda)}{W(f_-,f_+)(\lambda)} \right ]&=\mathrm{Im}\left [-b_{00}\frac{a_1^- a_1^+}{a_0^- a_0^+}+b_{01}\frac{a_1^- a_0^+}{a_0^- a_0^+}+
b_{10}\frac{a_0^- a_1^+}{a_0^- a_0^+}-b_{11} \right ]^{-1}(\lambda) \\
&= \frac{b_{01}+b_{10}}{b^2_{11}}\lambda+O(\lambda^{\alpha-1})
\end{align*}
and, analogously,
\begin{align*}
 \mathrm{Im}\left [\frac{a_0^- a_1^+(\lambda)}{W(f_-,f_+)(\lambda)} \right ]&=
\frac{1}{b_{11}}\lambda+O(\lambda^{\alpha-1}) \\
\mathrm{Im}\left [\frac{a_1^- a_0^+(\lambda)}{W(f_-,f_+)(\lambda)} \right ]&=
\frac{1}{b_{11}}\lambda+O(\lambda^{\alpha-1}) \\
\mathrm{Im}\left [\frac{a_1^- a_1^+(\lambda)}{W(f_-,f_+)(\lambda)} \right ]&=
O(\lambda^{\alpha-1}).
\end{align*}
The assertion for $\alpha=3$ follows by a similar computation.
\end{proof}

\section{Oscillatory integral estimates}
Based on Proposition \ref{prop_spectralmeasure} we are now able to finish the proof of Theorem \ref{thm_main}. Recall that we have to estimate terms of the form
$$ \lim_{N \to \infty}\int_\mathbb{R} \int_{1/N}^N \sin(t\lambda) \mathrm{Im}\left [\frac{f_-(x,\lambda)f_+(x',\lambda)}{W(f_-,f_+)(\lambda)} \right ]d\lambda f(x')dx'. $$
We do this by considering different regimes of $\lambda$, $x\lambda$ and $x'\lambda$ separately. To this end we introduce a smooth cut--off $\chi_\delta$ that satisfies $\chi_\delta(x)=1$ for $|x| \leq \delta$ and $\chi_\delta(x)=0$ for $|x| \geq 2\delta$ where $\delta>0$ is sufficiently small.
Since the estimates in this section are completely analogous to \cite{donninger} we only sketch the proofs and point out deviations from \cite{donninger}.

\subsection{Estimates for $\lambda$, $|x\lambda|$ and $|x'\lambda|$ small}
As always in the context of dispersive estimates, the most important contributions come from the small energy regime. 
We begin by showing that $u_j^\pm(x,\lambda)$ are even functions of $\lambda$.

\begin{lemma}
\label{lem_ujeven}
The functions $u_j^\pm(x,\lambda)$, $j=0,1$, constructed in Lemma \ref{lem_pertenergy}, satisfy the bounds
\begin{align*} 
|\partial_\lambda^{2k} u_0^\pm(x,\lambda)|&\leq C_k \langle x \rangle^{2k+1} \\
|\partial_\lambda^{2k+1} u_0^\pm(x,\lambda)|&\leq C_k \langle x \rangle^{2k+3}\lambda \\
|\partial_\lambda^{2k} u_1^\pm(x,\lambda)|&\leq C_k \langle x \rangle^{2k} \\
|\partial_\lambda^{2k+1} u_1^\pm(x,\lambda)|&\leq C_k \langle x \rangle^{2k+2}\lambda 
\end{align*}
for $k \in \mathbb{N}_0$, $\lambda>0$ small and $x \in [0,\lambda^{-1}]$ ($x \in [-\lambda^{-1},0]$, respectively).
\end{lemma}

\begin{proof}
For $k=0$ the claim follows directly from Lemma~\ref{lem_pertenergy}.
For $k\geq 1$ one proceeds by induction, cf.~Lemma 4.4 in \cite{donninger}. 
\end{proof}

In order to handle the small energy regime, the following general observation will be useful.

\begin{lemma}
\label{lem_osc}
Let $\omega: [0,\infty) \to \mathbb{R}$ be of the form
$$ \omega(\lambda)=\tilde{\omega}(\lambda)+O(\lambda^{\alpha-1}) $$
where $\alpha \geq 1$, $\tilde{\omega}$ is a smooth odd function and the $O$--term behaves like a symbol.
Moreover, assume that $\omega(\lambda)=0$ for all $\lambda \geq \lambda_0>0$.
Then we have
$$ \int_0^\infty \sin(t\lambda)\omega(\lambda)d\lambda \leq C(\omega)\langle t \rangle^{-\alpha} $$
for all $t \geq 0$. 
Furthermore, $C(\omega)$ can be estimated as
$$ C(\omega)\leq C\max \left \{\|\omega\|_{L^\infty(0,\infty)}, \|\omega^{([\alpha]-1)}\|_{L^\infty(0,\infty)}, \sup_{\lambda>0}|\lambda \omega^{([\alpha])}(\lambda)| \right \} $$
where $C>0$ is an absolute constant. 
\end{lemma}

\begin{proof}
It suffices to prove the lemma for $t \geq 1$. 
We integrate by parts $([\alpha]-1)$--times and obtain
$$ \left | \int_0^\infty \sin(t\lambda)\omega(\lambda)d\lambda \right |=\frac{1}{t^{[\alpha]-1}} \left |\int_0^\infty f(t\lambda)\omega^{([\alpha]-1)}(\lambda)d\lambda \right |$$
where $f$ stands for $\sin$ or $\cos$ depending on whether $[\alpha]$ is odd or even.
Note that all boundary terms vanish since $\tilde{\omega}$ is odd.
Let $\chi$ be a smooth cut--off supported in, say, $[0,2]$ with $\chi=1$ on $[0,1]$.
Then we have
$$ \left |\int_0^\infty f(t\lambda)\omega^{([\alpha]-1)}(\lambda)\chi(t\lambda)d\lambda \right |
=\frac{1}{t} \left | \int_0^\infty f(\eta)\underbrace{\omega^{([\alpha]-1)}(\eta/t)}_{\lesssim |\eta/t|^{\alpha-[\alpha]}}\chi(\eta)d\eta \right | \lesssim t^{[\alpha]-\alpha-1} $$
since $\alpha-[\alpha]>-1$.
Furthermore, an additional integration by parts yields
\begin{align} 
\label{eq_proofosc}
\left |\int_0^\infty f(t\lambda)\omega^{([\alpha]-1)}(\lambda)[1-\chi(t\lambda)]d\lambda \right | \leq &
\frac{1}{t} \left |\int_0^\infty f(t\lambda)\omega^{([\alpha])}(\lambda)[1-\chi(t\lambda)]d\lambda \right | \\
&+\left |\int_0^\infty f(t\lambda)\omega^{([\alpha]-1)}(\lambda)\chi'(t\lambda)d\lambda \right | \nonumber
\end{align}
where the boundary term vanishes thanks to the cut--off and the fact that $\omega$ is identically zero for large $\lambda$.
Now we estimate the two terms on the right--hand side of Eq.~\eqref{eq_proofosc} as
$$
\frac{1}{t} \left |\int_0^\infty f(t\lambda)\omega^{([\alpha])}(\lambda)[1-\chi(t\lambda)]d\lambda \right |
=\frac{1}{t^2} \left |\int_0^\infty f(\eta)\underbrace{\omega^{([\alpha])}(\eta/t)}_{\lesssim |\eta/t|^{\alpha-[\alpha]-1}}[1-\chi(\eta)]d\eta \right | \lesssim t^{[\alpha]-\alpha-1}
$$
since $\omega$ vanishes for large $\lambda$ and we stay away from the origin due to the cut--off.
For the second term we have
$$
\left |\int_0^\infty f(t\lambda)\omega^{([\alpha]-1)}(\lambda)\chi'(t\lambda)d\lambda \right |=\frac{1}{t}
\left |\int_0^\infty f(\eta)\underbrace{\omega^{([\alpha]-1)}(\eta/t)}_{\lesssim |\eta/t|^{\alpha-[\alpha]}}\chi'(\eta)d\eta \right | 
\lesssim t^{[\alpha]-\alpha-1}
$$
since $\mathrm{supp}\chi' \subset [1,2]$.
\end{proof}

\begin{lemma}
Let $\delta>0$ be sufficiently small. Then we have the bound
$$ \sup_{x,x' \in \mathbb{R}}\left |\int_0^\infty \sin(t\lambda)\mathrm{Im}\left [\frac{f_-(x,\lambda)f_+(x',\lambda)}{W(f_-,f_+)(\lambda)} \right ]\langle x \rangle^{-\alpha-1}\langle x' \rangle^{-\alpha-1}\chi_\delta(\lambda)\chi_\delta(x\lambda)\chi_\delta(x'\lambda)d\lambda \right |\lesssim \langle t \rangle^{-\alpha}. $$
for all $t \geq 0$.
\end{lemma}

\begin{proof}
In view of the representation Eq.~(\ref{eq_Jostrep}), the symmetry in $x$, $x'$ and by expressing $u_j^\pm(\cdot,\lambda)$ in terms of $u_0^\mp(\cdot,\lambda)$, $u_1^\mp(\cdot,\lambda)$, 
it suffices to show
$$ \sup_{x<0,x'>0}\left |\int_0^\infty \sin(t\lambda)\omega(x,x',\lambda)d\lambda \right |\lesssim \langle t \rangle^{-\alpha} $$
where
$$ \omega(x,x',\lambda)=\mathrm{Im}\left [\frac{a_j^-(\lambda)a_k^+(\lambda)}{W(f_-,f_+)(\lambda)} \right ]u_{j'}^-(x,\lambda)u_{k'}^+(x',\lambda)\langle x \rangle^{-\alpha-1}\langle x' \rangle^{-\alpha-1}\chi_\delta(\lambda)\chi_\delta(x\lambda)\chi_\delta(x'\lambda). $$
According to Proposition~\ref{prop_spectralmeasure} and Lemma~\ref{lem_ujeven}, $\omega$ decomposes as
$$ \omega(x,x',\lambda)=\tilde{\omega}(x,x',\lambda)+O(\langle x \rangle^{-\alpha}\langle x' \rangle^{-\alpha}\lambda^{\alpha-1}) $$
where $\tilde{\omega}(x,x',\lambda)$ satisfies 
\begin{align*}
|\partial_\lambda^{2\ell}\tilde{\omega}(x,x',\lambda)|&\leq C_\ell 
\langle x \rangle^{-\alpha+2\ell}\langle x' \rangle^{-\alpha+2\ell}\lambda \\
|\partial_\lambda^{2\ell+1}\tilde{\omega}(x,x',\lambda)|&\leq C_\ell 
\langle x \rangle^{-\alpha+2\ell}\langle x' \rangle^{-\alpha+2\ell}
\end{align*} 
for $\ell \in \mathbb{N}_0$ and the $O$--term behaves like a symbol with respect to differentiation in $\lambda$.
Furthermore, $\omega(x,x',\lambda)$ is compactly supported with respect to $\lambda$ and the claim follows by applying Lemma \ref{lem_osc}.
\end{proof}

\subsection{Estimates for the remaining cases}
\begin{lemma}
\label{lem_osc2}
 Let $\delta>0$ be sufficiently small. Then we have the estimates
$$ \sup_{x,x' \in \mathbb{R}}\left |\int_0^\infty \sin(t\lambda)\mathrm{Im}\left [\frac{f_-(x,\lambda)f_+(x',\lambda)}{W(f_-,f_+)(\lambda)} \right ]\langle x \rangle^{-\alpha-1}\langle x' \rangle^{-\alpha-1}\chi_\delta(\lambda)[1-\chi_\delta(x\lambda)][1-\chi_\delta(x'\lambda)]d\lambda \right |\lesssim \langle t \rangle^{-\alpha}. $$
\end{lemma}

\begin{proof}
The proof is similar to the corresponding case in~\cite{donninger}.
It is convenient to change the domain of integration from $\lambda>0$ to $\lambda \in \mathbb{R}$ which can be done easily since $\mathrm{Im}[G(x,x',\lambda)]$, the imaginary part of the Green's function, is odd in $\lambda$. This follows immediatly from the definition of the Jost functions.
As a consequence, we may replace 
$$\mathrm{Im}\left [\frac{f_-(x,\lambda)f_+(x',\lambda)}{W(f_-,f_+)(\lambda)} \right ]$$ in the oscillatory integral by 
$$ \frac{f_-(x,\lambda)f_+(x',\lambda)}{W(f_-,f_+)(\lambda)}$$ and change the domain of integration from the postive real axis to all of $\mathbb{R}$.

We distinguish different ranges of $x$, $x'$ and start with the case $x \leq 0$ and $x' \geq 0$.
By setting $m_\pm(x,\lambda)=e^{\mp i \lambda x}f_\pm(x,\lambda)$ we remove the oscillations and define
$$ \omega(x,x',\lambda):=\left [\frac{m_-(x,\lambda)m_+(x',\lambda)}{W(f_-,f_+)(\lambda)} \right ]\chi_\delta(\lambda)[1-\chi_\delta(x\lambda)][1-\chi_\delta(x'\lambda)]. $$
Thus, it suffices to estimate
$$ \sup_{x<0,x'>0}\left |\int_\mathbb{R} e^{i\lambda(\pm t-x+x')}\omega(x,x',\lambda)\langle x \rangle^{-\alpha-1}\langle x' \rangle^{-\alpha-1}d\lambda \right |$$
as $t \to \infty$.
According to Lemma \ref{lem_symbf+} and the symbol behavior of $W(f_-,f_+)(\lambda)$ we have the bound
$$ |\partial_\lambda^\ell \omega(x,x'\lambda)| \lesssim \langle x \rangle^{\ell}\langle x' \rangle^{\ell} $$
since $|\lambda^{-1}|\lesssim \langle x \rangle$ and hence, by $[\alpha]$--fold integration by parts, we obtain
$$ \sup_{x<0,x'>0}\left |\int_\mathbb{R} e^{i\lambda(\pm t-x+x')}\omega(x,x',\lambda)\langle x \rangle^{-\alpha-1}\langle x' \rangle^{-\alpha-1}d\lambda \right |\lesssim |\pm t-x+x'|^{-[\alpha]} \lesssim \langle t \rangle^{-\alpha} $$
provided that $|\pm t-x+x'|\geq \frac{1}{2}t$.
In the case $|\pm t-x+x'|\leq \frac{1}{2}t$ we have $\langle x \rangle^{-\alpha-1}\langle x' \rangle^{-\alpha-1}\lesssim \langle t \rangle^{-\alpha-1}$ which implies the same estimate.

In order to deal with terms that
involve $f_-(x,\lambda)$ for $x \geq 0$ or $f_+(x',\lambda)$ for $x' \leq 0$
we use reflection and transmission coefficients.
For $\lambda \not=0$, the functions $f_+(\cdot,\lambda)$ and
$\overline{f_+(\cdot,\lambda)}$ are linearly independent which shows that
there exist coefficients $a(\lambda)$ and $b(\lambda)$ such that
$f_-(x,\lambda)=a(\lambda)f_+(x,\lambda)+b(\lambda)\overline{f_+(x,\lambda)}
$.
This representation implies $|b(\lambda)|^2-|a(\lambda)|^2=1$ (cf.~
\cite{deift}) and thus,
$f_+(x,\lambda)=-\overline{a(\lambda)}f_-(x,\lambda)+b(\lambda)
\overline{f_-(x,\lambda)}$.
Furthermore, we have
$W(f_-,f_+)(\lambda)=b(\lambda)W(\overline{f_+},f_+)(\lambda)=2i\lambda b(\lambda)$ which
is equivalent to
\begin{equation}
\label{eq_b}
\frac{b(\lambda)}{W(f_-,f_+)(\lambda)}=\frac{1}{2i\lambda}.
\end{equation}
Similarly, we obtain
$ W(f_-,\overline{f_+})(\lambda)=-2i\lambda a(\lambda)$ and therefore,
$$ \frac{a(\lambda)}{W(f_-,f_+)(\lambda)}=-
\frac{W(f_-,\overline{f_+})(\lambda)}{2i\lambda
W(f_-,f_+)(\lambda)}. $$
However, from the representation of $W(f_-,f_+)(\lambda)$ given in the proof of Lemma \ref{prop_spectralmeasure} it follows that
$$ \frac{W(f_-,\overline{f_+})(\lambda)}{W(f_-,f_+)(\lambda)}=O_\mathbb{C}(1) $$
and hence,
\begin{equation}
 \label{eq_a}
 \frac{a(\lambda)}{W(f_-,f_+)(\lambda)}=O_\mathbb{C}(\lambda^{-1}) 
\end{equation}
where the $O_\mathbb{C}$--term behaves like a symbol.
This shows that, for $x \geq 0$ or $x' \leq 0$, one picks up an additional $\lambda^{-1}$ factor, see Eqs.~\eqref{eq_a} and \eqref{eq_b}.
However, negative powers of $\lambda$ can be controlled by positive powers of $\langle x' \rangle$ and one can proceed as above in the case $x\leq 0$ and $x'\geq 0$.
\end{proof}

Next, we consider the case $\lambda$, $|x\lambda|$ small and $|x'\lambda|$ large. By symmetry, this is equivalent to $\lambda$, $|x'\lambda|$ small and $|x\lambda|$ large.

\begin{lemma}
 Let $\delta>0$ be sufficiently small. Then we have the estimates
$$ \sup_{x,x' \in \mathbb{R}}\left |\int_0^\infty \sin(t\lambda)\mathrm{Im}\left [\frac{f_-(x,\lambda)f_+(x',\lambda)}{W(f_-,f_+)(\lambda)} \right ]\langle x \rangle^{-\alpha-1}\langle x' \rangle^{-\alpha-1}\chi_\delta(\lambda)\chi_\delta(x\lambda)[1-\chi_\delta(x'\lambda)]d\lambda \right |\lesssim \langle t \rangle^{-\alpha}. $$
\end{lemma}

\begin{proof}
 The proof is again completely analogous to the corresponding case in Section~8 of \cite{donninger}.
As before, we remove the oscillation from $f_+(x',\lambda)$ and use the bounds from Lemma~\ref{lem_symbf+}.
However, for $f_-(x,\lambda)$ we use the representation Eq.~(\ref{eq_Jostrep}) since $|x\lambda|$ is small.
Again, we distinguish between $x'\geq 0$ and $x'\leq 0$ where the latter case is more difficult since one picks up a $|\lambda|^{-1}$ factor, see Eqs.~\eqref{eq_a} and \eqref{eq_b}.
For instance, we have to deal with the term
$$ \omega(x,x',\lambda):=\frac{a_0^-(\lambda)u_1^-(x,\lambda)\overline{a(\lambda)}m_-(x',\lambda)}{W(f_-,f_+)(\lambda)}
\chi_\delta(\lambda)\chi_\delta(x\lambda)[1-\chi_\delta(x'\lambda)] $$
and we need to estimate 
$$ \sup_{x<0,x'<0}\left |\int_\mathbb{R}e^{i\lambda(\pm t-x')}\omega(x,x',\lambda)\langle x \rangle^{-\alpha-1}\langle x' \rangle^{-\alpha-1}d\lambda \right |. $$
In the domain considered, $\omega$ satisfies the bounds
$$ |\partial_\lambda^\ell \omega(x,x',\lambda)|\leq C_\ell \langle x \rangle^{\ell+1} \langle x' \rangle^{\ell+1} $$
for all $\ell \in \mathbb{N}_0$
since one can trade $|\lambda|^{-1}$ for $\langle x' \rangle$.
The bound
$$ \sup_{x<0,x'<0}\left |\int_\mathbb{R}e^{i\lambda(\pm t-x')}\omega(x,x',\lambda)\langle x \rangle^{-\alpha-1}\langle x' \rangle^{-\alpha-1}d\lambda \right | \lesssim \langle t \rangle^{-\alpha} $$
now follows by appropriate integration by parts, see the proofs of Lemmas \ref{lem_osc} and \ref{lem_osc2}.
The remaining cases are completely analogous.
\end{proof}

Finally, we have to deal with the large energy contributions.
Clearly, those have nothing to do with the characteristic behavior of the spectral measure at zero energy and are more or less independent of the decay properties of the potential.
Therefore, we can directly transfer the results from Section~9 in~\cite{donninger} to our problem and conclude the proof of Theorem \ref{thm_main}.

\bibliography{cubic}

\begin{thebibliography}{10}

\bibitem{bizontail}
Piotr Bizo\'{n}, Tadeusz Chmaj, and Andrzej Rostworowski.
\newblock Anomalously small wave tails in higher dimensions.
\newblock {\em Physical Review D (Particles, Fields, Gravitation, and
  Cosmology)}, 76(12):124035, 2007.

\bibitem{ching}
E.~S.~C. Ching, P.~T. Leung, W.~M. Suen, and K.~Young.
\newblock Wave propagation in gravitational systems: Late time behavior.
\newblock {\em Physical Review D (Particles, Fields, Gravitation, and
  Cosmology)}, 52(4):2118--2132, 1995.

\bibitem{D-Rod}
Mihalis Dafermos and Igor Rodnianski.
\newblock A proof of {P}rice's law for the collapse of a self-gravitating
  scalar field.
\newblock {\em Invent. Math.}, 162(2):381--457, 2005.

\bibitem{deift}
P.~Deift and E.~Trubowitz.
\newblock Inverse scattering on the line.
\newblock {\em Comm. Pure Appl. Math.}, 32(2):121--251, 1979.

\bibitem{dss}
Roland Donninger, Wilhelm Schlag, and Avy Soffer.
\newblock On pointwise decay of linear waves on a {S}chwarzschild black hole
  background.
\newblock {\em Preprint}, 2009.

\bibitem{donninger}
Roland Donninger, Wilhelm Schlag, and Avy Soffer.
\newblock A proof of {P}rice's {L}aw on {S}chwarzschild black hole manifolds
  for all angular momenta.
\newblock {\em Preprint arXiv:0908.4292}, 2009.

\bibitem{Kron2}
Johann Kronthaler.
\newblock Decay rates for spherical scalar waves in the {S}chwarzschild
  geometry.
\newblock {\em Preprint arXiv:0709.3703}, 2007.

\bibitem{Price1}
Richard~H. Price.
\newblock Nonspherical perturbations of relativistic gravitational collapse.
  {I}. {S}calar and gravitational perturbations.
\newblock {\em Phys. Rev. D (3)}, 5:2419--2438, 1972.

\bibitem{Price2}
Richard~H. Price.
\newblock Nonspherical perturbations of relativistic gravitational collapse.
  {II}. {I}nteger-spin, zero-rest-mass fields.
\newblock {\em Phys. Rev. D (3)}, 5:2439--2454, 1972.

\bibitem{schlag3}
W.~Schlag.
\newblock Dispersive estimates for {S}chr\"odinger operators: a survey.
\newblock In {\em Mathematical aspects of nonlinear dispersive equations},
  volume 163 of {\em Ann. of Math. Stud.}, pages 255--285. Princeton Univ.
  Press, Princeton, NJ, 2007.

\bibitem{schlag1}
Wilhelm Schlag, Avy Soffer, and Wolfgang Staubach.
\newblock Decay for the wave and {S}chr\"odinger evolutions on manifolds with
  conical ends, {P}art {I}.
\newblock {\em Preprint arXiv:0801.1999, will appear in Trans. Amer. Math.
  Soc.}, 2008.

\bibitem{schlag2}
Wilhelm Schlag, Avy Soffer, and Wolfgang Staubach.
\newblock Decay for the wave and {S}chr\"odinger evolutions on manifolds with
  conical ends, {P}art {II}.
\newblock {\em Preprint arXiv:0801.2001, will appear in Trans. Amer. Math.
  Soc.}, 2008.

\bibitem{tataru}
Daniel Tataru.
\newblock Local decay of waves on asymptotically flat stationary space-times.
\newblock {\em Preprint arXiv:0910.5290}, 2009.

\bibitem{teschl}
Gerald Teschl.
\newblock {\em Mathematical methods in quantum mechanics}, volume~99 of {\em
  Graduate Studies in Mathematics}.
\newblock American Mathematical Society, Providence, RI, 2009.
\newblock With applications to Schr{\"o}dinger operators.

\end{thebibliography}
\bibliographystyle{plain}

\end{document}